\documentclass[12pt]{article}
\usepackage{amsfonts}
\usepackage{mathrsfs}
\usepackage{amsmath,amssymb}
\baselineskip 5.5pt \lineskip 5.5pt \numberwithin{equation}{section}

\def\QED{\hfill$\Box$\par}

\def\DD{\mathcal {D}}

\def\TT{\mathcal {T}}

\def\VV{\mathcal {V}}
\def\UU{\mathcal {U}}
\def\CL{\mathcal {C}}

\def\LL{\mathcal {L}}

\def\cl{\centerline}
\def\DD{{\cal D}}

\def\rar{\longrightarrow}

\def\vs{\vspace*}

\def\WW{{\cal W}}
\def\AA{{\cal A}}

\def\C{\mathbb{C}}

\def\Z{\mathbb{Z}}

\newtheorem{theo}{Theorem}[section]
\newtheorem{coro}[theo]{Corollary}
\newtheorem{lemm}[theo]{Lemma}

\newtheorem{defi}[theo]{Definition}

\begin{document}

\cl{\bf Quantum group structure of the $q$-deformed $W$ algebra
$\WW_q$\footnote {Supported by NSF grants 10471091, 10671027 of
China, ``One Hundred Talents Program'' from University of Science
and Technology of China.\\[2pt] \indent Corresponding E-mail:
sd\_huanxia@163.com}} \vs{6pt}

\cl{Huanxia Fa$^{\dag)}$, Junbo Li$^{\dag)}$, Yongsheng
Cheng$^{\ddag)}$} \cl{\small $^{\dag)}$Department of Mathematics,
Changshu Institute of Technology, Changshu 215500, China} \cl{\small
$^{\ddag)}$College of Mathematics and Information Science, Henan
University, Kaifeng 475001, China} \cl{\small E-mail:
sd\_huanxia@163.com, sd\_junbo@163.com, yscheng@henu.edu.cn}
\vs{6pt}

\noindent{{\bf Abstract.} In this paper the $q$-deformed $W$ algebra
$\WW_q$ is constructed, whose nontrivial quantum group structure is
presented. \vs{5pt}

}

\parskip .001 truein\baselineskip 8pt \lineskip 8pt

\vs{6pt}
\par
\cl{\bf\S1. \ Introduction}
\setcounter{section}{1}\setcounter{theo}{0}\setcounter{equation}{0}
\par

The deformation theory acts important roles in many fields such as
mathematics and physics, which is closely related to quantum groups,
originally introduced by Drinfeld in \cite{D1}. From the day that
the conception of the quantum groups was born, there appear many
papers on this relatively new object, so does the deformation theory
(cf. \cite{HN,AS,DG,FR}, \cite{HLS}--\cite{SW}). The quantum group
structure on the $q$-deformed Virasoro algebra and $q$-deformed
Kac-Moody algebra had been investigated by many authors (cf.
\cite{HN}, \cite{HLS}--\cite{SW}), and some interesting results were
presented therein. In particular, the structure and representations
of $q$-Virasoro algebra were intensively investigated in \cite{HN}.
In \cite{CS}, $q$-deformation of the twisted Heisenberg-Virasoro
algebra with central extension was constructed, which admitted a
nontrivial Hopf structure.

Now let's introduce the object algebra concerned with in the present
paper. The algebra $W$-algebra $W(2,2)$, denoted by $\WW$ for
convenience and introduced by Zhang and Dong in \cite{ZD}, is an
infinite-dimensional Lie algebra, possessing a $\C$-basis
$\{\,L_n,\,W_n,\,\mathcal {C}\,|\,n\in \Z\,\}$ and admitting the
following Lie brackets (\,other components vanishing):
\begin{eqnarray}
&&[L_m,L_n]=(m-n)L_{m+n}+\frac{m^3-m}{12}\delta_{m+n,0}\CL,\label{BLie01}\\
&&[L_m,W_n]=(m-n)W_{m+n}+\frac{m^3-m}{12}\delta_{m+n,0}\CL.\label{BLie02}
\end{eqnarray}
There appeared some papers investigating the structures and
representations on such $W$ algebra recently. In \cite{ZD}, Zhang
and Dong produced a new class of irrational vertex operator algebras
by studying its highest weight modules, while \cite{LGZ} and
\cite{LZ} classified its irreducible weight modules and
indecomposable modules and \cite{GJP} determined its derivations,
central extensions and automorphisms. Afterwards, the Lie bialgebra
structures on $\WW$ (centerless form) were proved to triangular
coboundary in \cite{LS1}, which were quantized in \cite{LS2}.
However, the existence of a $q$-deformation of the $W$-algebra
$W(2,2)$ and its quantum group structure is still an open problem,
which may be interesting to physicists. That is what our paper shall
focus on. In other words, we shall construct a $q$-deformation of
the $W$ algebra, i.e., $\WW_q$, which admits a nontrivial Hopf
structure. The Harish-Chandra modules, Verma modules and also
Unitary representations for the $q$-deformed $W$-algebra $\WW_q$
have been investigated and shall be presented in a series of papers
(c.f. \cite{LS3}--\cite{LS5}).

Let's formulate our main results below. The following definition can
be found in many references (e.g. \cite{HN}, \cite{CS}, \cite{HLS}).
\begin{defi}\label{def1}
A vector space $\VV$ over $\C$, with an bilinear operation
$\VV\times\VV\rar\VV$, denoted $(x,y)\rar[x,y]_q$ and called the
$q$-bracket or $q$-commutator of $x$ and $y$, and meanwhile with an
endomorphism of $\VV$, denoted $f_q$, is called a $q$-deformed Lie
algebra over $\C$ if the following axioms are satisfied:\vs{-8pt}
\begin{eqnarray}
&&[u,v]_q=-[v,u]_q,\label{Lde1} \\
&&\big[f_q(u),[v,w]_q\big]_q+\big[f_q(w),[u,v]_q\big]_q+\big[f_q(v),[w,u]_q\big]_q=0,\label{Lde2}
\end{eqnarray}
for any $u,v,w\in\VV$.
\end{defi}

As the usual definition of $2$-cocycle, we also can introduce the
corresponding one of $q$-deformed 2-cocycle on the centerless
$q$-deformed Lie algebra $\VV$ defined in Definition \ref{def1}.
\begin{defi}\label{def2}
A bilinear $\C$-value function $\psi_q:\VV\times\VV\rar \C$ is
called q-deformed 2-cocycle on $\VV$ if the following conditions are
satisfied
\begin{eqnarray}
&&\psi_q(u,v)=-\psi_q(v,u),\label{Lcoy1}\\
&&\psi_q\big(f_q(u),[v,w]_q\big)+\psi_q\big(f_q(w),[u,v]_q\big)
+\psi_q\big(f_q(v),[w,u]_q\big)=0,\label{Lcoy2}
\end{eqnarray}
for any $u,v,w\in\VV$.
\end{defi}

Denote by $\mathscr{C}_q^2(\VV,\C)$ the vector space of $q$-deformed
2-cocycles on $\VV$. For any linear $\C$-value function
$\chi_q:\VV\rar\C$, the 2-cocycle $\psi_{\chi_q}$ defined by
\begin{equation}\label{la1.8}
\psi_{\chi_q}(u, v)=\chi_q([u,v]_q),\ \ \forall\,\,u,v\in\VV,
\end{equation}
\noindent is called 2-coboundary on $\VV$. Denote by $\mathscr{B}^2
(\VV,\C)$ the vector space of 2-coboundaries on $\VV$. The quotient
space $\mathscr{H}^2(\VV,\C):=\mathscr{C}^2
(\VV,\C)/\mathscr{B}^2(\VV,\C)$ is called the second cohomlogy group
of $\VV$.

\begin{theo}\label{mainth}
The algebra $\UU_q$ is a noncommutative but cocommutative Hopf
algebra under the comultiplication $\Delta$, the counity $\epsilon$
and the antipode $\mathcal {S}$ defined by
(\ref{dhq1})--(\ref{dhq3}).
\end{theo}

\cl{\bf\S2. \ Proof of the main result
}\setcounter{section}{2}\setcounter{theo}{0}\setcounter{equation}{0}

Firstly, we shall construct a $q$-deformation of $\WW$, denoted
$\WW_q$, by using some techniques developed in \cite{HN,CS,HLS}. In
fact, the Witt algebra can be recognized as the Lie algebra of
derivations on $\mathbb{C}[t^{\pm1}]$, i.e., the Lie algebra of its
linear operators $\Omega$ satisfying
$$\Omega(xy)=\Omega(x)y+x\Omega(y),$$ whose Lie
bracket also can be obtained by simple computations. Fix some
generic $q\in\mathbb{C}^*$, and
$\delta\in\mbox{End}(\mathbb{C}[t^{\pm1}])$ such that
$\delta(t)=qt$. Define a $q$-derivation $\DD$ as
\begin{eqnarray*}
\DD(f(t))=-{(q-1)}^{-1}{(Id-\delta)}f(t), \ \forall\,\,f(t)\in
\mathbb{C}[t^{\pm1}].
\end{eqnarray*}
It is easy to see that $\delta(t^n)=q^nt^n$ and
$\DD(t^n)=\frac{q^nt^n-t^n}{q-1}=[n]_qt^n$, where
$[n]_q=\frac{q^n-1}{q-1}$ for some $n\in\mathbb{Z}$. The following
way of defining $q$-deformed Virasoro algebra can be found in many
references (e.g., \cite{HN,HLS}), on which our construction is
based\vs{-8pt}
\begin{eqnarray}\label{LBqLL}
[L_m,L_n]_q=([m]_q-[n]_q)L_{m+n}+
\frac{q^{-m}[m-1]_q[m]_q[m+1]_q}{6(1+q^m)}\delta_{m,-n}\CL.
\end{eqnarray}
\begin{defi}\label{defvico}
The 2-cocycle on the $q$-deformed Virasoro algebra given in
(\ref{LBqLL}) is called the $q$-deformed Virasoro 2-cocycle.
\end{defi}
Combining the structures of the algebra $\WW$ listed in
(\ref{BLie01})--(\ref{BLie02}) and the $q$-deformed Virasoro Lie
algebras given in (\ref{LBqLL}), we introduce the centerless
$q$-deformed $W$ algebra $\mathscr{W}_q$, which possesses a
$\C$-basis $\{L_m,W_m\,|\,m\in\Z\}$ with the following relations
\begin{eqnarray}\label{LBLWmn}
[L_m,L_n]_{q}=([m]_q-[n]_q)L_{m+n},\ \
[L_m,W_n]_{q}=([m]_q-[n]_q)W_{m+n},\ \ [W_m,W_n]_{q}=0.
\end{eqnarray}
Observing (\ref{Lde2}), (\ref{Lcoy2}), (\ref{LBqLL}) and
(\ref{LBLWmn}), one can take
\begin{eqnarray}\label{eqde1}
f_q(L_m)=(q^m+1)L_m,\, \ f_q(W_m)=(q^m+1)W_m, \
\,\forall\,\,m\in\mathbb{Z},
\end{eqnarray}
where $f_q$ is that defined in Definition \ref{def1}. By simple
computations, one can see that the algebra $\mathscr{W}_q$ defined
by (\ref{LBLWmn}) with the $f_q$ defined by (\ref{eqde1}) is indeed
a $q$-deformed Lie algebra.

Using (\ref{LBqLL}), in order to obtain the $q$-deformed algebra
$\WW_q$, we have to determine the $q$-deformed 2-cocycle
$\psi_q(L_m,W_n)$ determined by the following identity
\begin{eqnarray}\label{dq1}
[L_m,W_n]_{q}=([m]_q-[n]_q)W_{m+n}+\psi_q(L_m,W_n)\,\CL.
\end{eqnarray}

Using (\ref{eqde1}) and respectively, replacing $(u,v)$ by
$(L_i,W_j)\,(\forall\,i,j\in\Z)$ in (\ref{Lcoy1}) and the triple
$(u,v,w)$ by $(L_i,L_j,W_k)\,(\forall\,i,j,k\in\Z)$ in
(\ref{Lcoy2}), one has
\begin{eqnarray}
&&\psi_q(L_i,W_j)=-\psi_q(W_j,L_i),\label{Lcoy1-1}\\
&&(q^i+1)([j]_q-[k]_q)\psi_q\big(L_i,W_{j+k}\big)\nonumber\\
&&=(q^k+1)([i]_q-[j]_q)\psi_q\big(L_{i+j},W_k\big)
+(q^j+1)([i]_q-[k]_q)\psi_q\big(L_j,W_{k+i}\big).\label{Lcoy2-1}
\end{eqnarray}
Let $i=0$ in (\ref{Lcoy2-1}), one has
\begin{eqnarray*}
&&(q^j-q^k)\psi_q\big(L_0,W_{j+k}\big)
=(q^{j+k}-1)\psi_q\big(L_{j},W_k\big),
\end{eqnarray*}
which together with our assumption on $q$, forces
\begin{eqnarray}
&&\psi_q\big(L_0,W_{0}\big)=0.\label{LcoyLW00}
\end{eqnarray}
According to the second bracket in (\ref{LBLWmn}), we can write
\begin{eqnarray*}
&&L_0=(1+q^{-1})[L_{1},L_{-1}]_q,\ \ W_0=(1+q^{-1})[L_{1},W_{-1}]_q,\\
&&L_m=\big([m]_q\big)^{-1}[L_0,L_m]_q,\ \
W_m=\big([m]_q\big)^{-1}[L_0,W_m]_q\ \ {\rm if}\ \,m\in\Z^*.
\end{eqnarray*}
Define a $\C$-linear function $\chi_q:\mathscr{W}_q\rightarrow\C$ as
follows
\begin{eqnarray*}
&&\chi_q(L_0)=(1+q^{-1})\psi_q(L_{1},L_{-1}),\ \ \chi_q(W_0)=(1+q^{-1})\psi_q(L_{1},W_{-1}),\\
&&\chi_q(L_m)=\big([m]_q\big)^{-1}\psi_q(L_0,L_m),\ \
\chi_q(W_m)=\big([m]_q\big)^{-1}\psi_q(L_0,W_m)\ \ {\rm if}\
\,m\in\Z^*.
\end{eqnarray*}
Let $\varphi_q=\psi_q-\psi_{\chi_q}$ where $\psi_{\chi_q}$ is
defined in (\ref{la1.8}). One has
\begin{eqnarray}
&&\varphi_q(L_{1},L_{-1})=\varphi_q(L_{1},W_{-1})=\varphi_q(L_0,L_m)=\varphi_q(L_0,W_m)=0\
\ {\rm if}\ \,m\in\Z^*.\label{reBa2}
\end{eqnarray}
Denote by $\mathfrak{W}_q$ the $q$-deformed Witt subalgebra of
$\mathscr{W}_q$ spanned by $\{L_m\,|\,m\in\Z\}$. The by simple
discussion or cite the result given in \cite{HN,HLS}, one can
suppose that $\varphi_q|_{\mathfrak{W}_q}$ is exactly the
$q$-deformed Virasoro 2-cocycle (up to a constant factor).

Recalling (\ref{LcoyLW00}) and (\ref{reBa2}), one can deduce
$\varphi_q(L_m,W_n)=0$ if $m+n\neq 0$. Thus, the left components we
have to compute are
\begin{eqnarray}
&&\varphi_q(L_{m},W_{-m}),\ \,\ \,\forall\,\,m\in\Z^*.\label{reBa2}
\end{eqnarray}
By employing the same techniques developed in \cite{HN,HLS}, we
obtain (up to a constant factor)
\begin{eqnarray}
&&\varphi_q(L_m,W_{-m})=\frac{q^{-m}[m-1]_q[m]_q[m+1]_q}{6(1+q^m)},\
\,\ \,\forall\,\,m\in\Z^*.\label{reBa2-2}
\end{eqnarray}
Then we have
\begin{eqnarray}\label{LBqLW}
[L_m,W_n]_q=([m]_q-[n]_q)W_{m+n}+
\frac{q^{-m}[m-1]_q[m]_q[m+1]_q}{6(1+q^m)}\delta_{m,-n}\CL.
\end{eqnarray}
Now we can safely present the following lemma.
\begin{lemm}\label{lem2.2}
The algebra $\WW_q$ with a $\C$-basis
$\{L_m,W_m,\CL\,|\,m\in\mathbb{Z}\}$ satisfying the following
relations (while other components vanishing) is a $q$-deformation of
the algebra $\WW$.
\begin{eqnarray}\label{LBwa1}
[L_m,L_n]_q=q^{m}L_{m}L_{n}-q^{n}L_{n}L_{m},\ \ \
[L_m,W_n]_q=q^{m}L_{m}W_{n}-q^{n}W_{n}L_{m},
\end{eqnarray}
where the $q$-deformed brackets are respectively given in
(\ref{LBqLL}) and (\ref{LBqLW}).
\end{lemm}

Next we shall proceed with our construction of the Hopf algebra
structure based on the $q$-deformed algebra $\WW_q$ given in Lemma
\ref{lem2.2}\,. Firstly, for convenience to express, we shall recall
the definition of a Hopf algebra, which can be found in many books
and also references.
\begin{defi}\label{deh}
A tuple $(\AA,\nabla,\varepsilon,\Delta,\epsilon,\,\mathcal {S})$,
$\AA$ being a $\mathbb{C}$-vector space, $\nabla:\AA\otimes
\AA\longrightarrow \AA$ a multiplication map,
$\varepsilon:\mathbb{C}\longrightarrow\AA$ a unit map,
$\Delta:\AA\longrightarrow \AA\otimes \AA$ a comultiplication map,
$\epsilon:\AA\longrightarrow \mathbb{C}$ a counit map, $\mathcal
{S}:\,\AA\longrightarrow \AA$ an antipode map, is called a Hopf
algebra over $\mathbb{C}$ if the following axioms are satisfied\\
(1) the map $\nabla$ gives an associative algebra structure
on $\AA$ with the unit $\varepsilon(1)$,\\
(2) $\Delta$ and $\epsilon$ give a coassociative coalgebra structure
on $\AA$,
\begin{eqnarray}
&&(1\otimes \Delta) \Delta(x)=(\Delta\otimes1)\Delta(x),\ \ \
(1\otimes
\epsilon)\Delta(x)=(\epsilon\otimes1)\Delta(x),\label{coaa}
\end{eqnarray}
(3) both $\Delta$ and $\epsilon$ are algebra homomorphisms,\\
(4) $\mathcal {S}$ is an automorphism with the following relations
\begin{eqnarray}\label{abc1}
\nabla(1\otimes\mathcal {S})\Delta(x)=\nabla(\mathcal
{S}\otimes1)\Delta(x)=\varepsilon(\epsilon(x)).
\end{eqnarray}
\end{defi}
We say the Hopf algebra $\AA$ is cocommutative if
$\Delta=\Delta^{op}$. A vector space $\LL$ over $\mathbb{C}$, is
called a bialgebra if it admits the maps
$\nabla,\varepsilon,\Delta,\,\epsilon$ with the axioms $(1)$--$(3)$
given in Definition \ref{deh}.

Denote $\UU_q$ to be the $q$-deformed enveloping algebra of $\WW_q$.
Then $\UU_q$ allows the Hopf algebra structure given below
\begin{eqnarray}
&&\epsilon(L_m)=\epsilon(W_m)=\epsilon(\CL)=0,\ \ \
\Delta(\CL)=\CL\otimes 1+1\otimes \CL,\label{dhq1}\\
&&\Delta(L_m)=L_m\otimes\TT^m+\TT^m\otimes L_m,\
\Delta(W_m)=W_m\otimes
\TT^m+\TT^m\otimes W_m,\label{dhq2}\\
&&\mathcal {S}(L_m)=-\TT^{-m}L_m\TT^{-m},\ \ \ \mathcal
{S}(W_m)=-\TT^{-m}W_m\TT^{-m},\ \ \ \mathcal
{S}(\CL)=-\CL,\label{dhq3}
\end{eqnarray}
where the operators $\{\TT,\TT^{-1}\}$ are given by
\begin{eqnarray}
&& \Delta(\TT)=\TT\otimes\TT, \ \ \epsilon(\TT)=1,\ \ \mathcal
{S}(\TT)=\TT^{-1}.\label{dhq5}
\end{eqnarray}
The following relations also can be obtained by simple computations:
\begin{eqnarray*}
&&\!\!\!\!\!\!\!\!\!\!\!\!\TT^mL_n=q^{-(n+1)m}L_n\TT^m,\ \
\TT^mW_n=q^{-(n+1)m}W_n\TT^m,\\
&&\!\!\!\!\!\!\!\!\!\!\!\!\TT^mL_n=q^{-(n+1)m}L_n\TT^m,\ \
\TT^mW_n=q^{-(n+1)m}W_n\TT^m,\\
&&\!\!\!\!\!\!\!\!\!\!\!\!\TT\TT^{-1}=\TT\TT^{-1}=1,\ \
q^m\TT^m\CL=\CL\TT^m,\ \ q^m\TT^m\CL=\CL\TT^m.
\end{eqnarray*}

{\bf Proof of Theorem \ref{mainth}}\ \ We shall follow some
techniques developed in \cite{H}. It is not difficult to see that
the coassociativity and cocommutative of $\Delta$ hold in $\UU_q$
and, $\epsilon$ is an algebra homomorphism, also
$(1\otimes\epsilon)\Delta=(\epsilon\otimes1)\Delta=1$. Firstly, We
shall ensure that $\Delta$ is an algebra homomorphism while
$\mathcal {S}$ is an algebra anti-homomorphism of $\UU_q$. Using the
relations obtained above, we can present the following computations:
\begin{eqnarray*}
&&q^{m}\Delta(L_{m})\Delta(W_{n})-q^{n}\Delta(W_{n})\Delta(L_{m})\\
&&=\big(q^{m}L_{m}W_{n}-q^{n}W_{n}L_{m}\big)\otimes
\TT^{m+n}+\TT^{m+n}\otimes
\big(q^{m}L_{m}W_{n}-q^{n}W_{n}L_{m}\big)\\
&&=[L_m,W_n]_q\otimes \TT^{m+n}+\TT^{m+n}\otimes [L_m,W_n]_q\\
&&=([m]_q-[n]_q)\Delta(W_{m+n})
+\frac{q^{-m}[m-1]_q[m]_q[m+1]_q}{6(1+q^m)}\delta_{m,-n}\Delta(\CL).
\end{eqnarray*}
Other formulate also can be proved to be preserved by the map
$\Delta$, which together implies that $\Delta$ is an algebra
homomorphism. Thus, $\UU_q$ indeed a bialgebra. We also have the
following computations:
\begin{eqnarray*}
\mathcal {S}(L_{m}W_{n})=\mathcal {S}(W_{n})\mathcal
{S}(L_{m})=\TT^{-n}W_n\TT^{-n}\TT^{-m}L_m\TT^{-m}
=q^{n-m}\TT^{-m-n}L_nL_m\TT^{-m-n},
\end{eqnarray*} which further gives
\begin{eqnarray*}
&&q^{m}\mathcal {S}(L_{m}W_{n})-q^{n}\mathcal {S}(W_{n}L_{m})\\
&&=q^{n}\TT^{-m-n}W_nL_m\TT^{-m-n}-q^{m}\TT^{-m-n}L_mW_n\TT^{-m-n}\\
&&=-\TT^{-m-n}(q^{m}L_mW_n-q^{n}W_nL_m)\TT^{-m-n}\\
&&=-\TT^{-m-n}[L_m,W_n]_q\TT^{-m-n}\\
&&=-([m]_q-[n]_q)\mathcal {S}(W_{m+n})
+\frac{q^{-m}[m-1]_q[m]_q[m+1]_q}{6(1+q^m)}\delta_{m,-n}\mathcal
{S}(\CL),
\end{eqnarray*}
and which actually implies the fact that $\mathcal {S}$ preserves
the second identity of (\ref{LBwa1}). Other formulate also can be
proved to be preserved by the antipode map $\mathcal {S}$. Thus,
$\UU_q$ admits the referred Hopf algebra structure.\QED

Before ending this short note, employing the main techniques
developed in \cite{H}, one can easily obtain the following
corresponding corollary.\vs{-6pt}
\begin{coro}
\label{coro} As vector spaces,
\begin{equation} \label{theo31}
\UU_q\cong\C[\TT,\TT^{-1}]\otimes_\mathbb{C}\UU(\WW_q),
\end{equation}
where $\UU(\WW_q)$ is the universal enveloping algebra of $\WW_q$
generated by $\{L_m,W_m,\CL\,|\,m\in\mathbb{Z}\}$ with the relations
presented in (\ref{LBwa1}).
\end{coro}

\end{document}